\documentclass[12pt]{amsart}
\input {epsf}
\usepackage{latexsym, amsbsy, amsmath, amsfonts, amssymb, amsthm, amscd, amscd}
\usepackage[all]{xy}
\newcommand{\Z}{\mathbb Z}

\newtheorem{Theorem}{Theorem}[section]

\newtheorem{Proposition}{Proposition}[section]

\allowdisplaybreaks[2]
\numberwithin{equation}{section}
\numberwithin{figure}{section}
\title{On the inverse mapping class monoids }
\author[Karoui]{R.~Karoui}
\author[Vershinin]{V.~V.~Vershinin}%$^{*}$}
\address{D\'epartement des Sciences Math\'ematiques,
                                     Universit\'e Montpellier II,
Place Eug\`ene Bataillon,
34095 Montpellier cedex 5, France}
\email{Rym.Karoui@math.univ-montp2.fr}
\address{D\'epartement des Sciences Math\'ematiques,
                                     Universit\'e Montpellier II,
Place Eug\`ene Bataillon,
34095 Montpellier cedex 5, France}
\email{ vershini@math.univ-montp2.fr}
\address{Sobolev Institute of Mathematics, Novosibirsk 630090,
Russia } 
\email{ versh@math.nsc.ru}

\begin{document}
\begin{abstract}
Braid groups and mapping class groups have many features in common.
Similarly to the notion of inverse braid monoid 
inverse mapping class monoid is defined. It concerns surfaces with punctures,
but among given $n$ punctures several can be omitted.
This corresponds to  braids where the number of strings
is not fixed. 
In the paper we give the analogue of 
the Dehn-Nilsen-Baer theorem, propose a
presentation of the inverse mapping class monoid for 
a punctured sphere and study the word problem.
This shows that certain  properties and objects based on 
mapping class
groups may be extended to the inverse mapping class monoids.
 We also give an analogues of Artin presentation
with two generators.
\end{abstract}
\maketitle
\tableofcontents

\section{Introduction}

Mapping class group is an  important object in Topology,
Complex Analysis, Algebraic Geometry and other domains. 
It is a lucky case when the method of Algebraic Topology works
perfectly well, the application of the functor of fundamental group
completely solves the topological problem: group of isotopy classes of
homeomorpisms
is described in terms of automorphisms of the fundamental group
of the corresponding surface, as states  the Dehn-Nilsen-Baer theorem, 
see \cite{Iv}, for example. 

Let  $S_{g,b,n}$ be an oriented surface of the genus $g$ with $b$ boundary components and
a set $Q_n$ of $n$ fixed points. Consider the group 
$\operatorname{Homeo}(S_{g,b,n})$ 
of orientation preserving
self-homeomorphisms of $S_{g,b,n}$ which fix pointwise the boundary (if
it exists) and map the set $Q_n$ into itself. 
Orientation reversing homeomorphisms also possible to consider, see 
\cite{DF}, for example,
but for the simplicity of exposition we restrict ourselves to orientation
preserving case. Let 
$\operatorname{Homeo}^0(S_{g,b,n})$ be the normal subgroup of 
self-homeomorphisms of $S_{g,b,n}$ which are isotopic to identity. 
Then the {\it mapping class group} $\mathcal {M}_{g,b,n}$ is defined as a factor group

\begin{equation*}
\mathcal {M}_{g,b,n} = \operatorname{Homeo}(S_{g,b,n})/
\operatorname{Homeo}^0(S_{g,b,n})
\end{equation*}

%Mapping class 
These groups are connected closely with braid groups.
In \cite{Mag1} W.~Magnus interpreted the braid group as the mapping class
group of a punctured disc with the fixed boundary.
Braid groups have a variety of generalizations, see \cite{Ve10}, for example.
One of generalizations is the {\it inverse braid monoid }  $IB_n$  constructed by D.~Easdown and T.~G.~Lavers \cite{EL}. 

The notion of {\it inverse semigroup} was introduced by V.~V.~Wagner in 1952 \cite{Wag}.
By definition it means that for any element $a$ of a semigroup (monoid)
$M$ there exists a unique element $b$ (which is called {\it inverse}) such that 
\begin{equation*}
a = aba
\label{eq:reg_v_n}
\end{equation*}
 and 
\begin{equation*}
b = bab.
\label{eq:inv}
\end{equation*}

The typical example of an inverse monoid is a monoid of partial (defined on a subset)
injections of a set. For a finite set this gives us the notion of a
{\it symmetric inverse monoid } $I_n$ which generalizes and includes the classical
symmetric group $\Sigma_n$. A presentation of symmetric inverse monoid was
obtained by L.~M.~Popova \cite{Po}, see also formulas 
(\ref{eq:brelations}), (\ref{eq:invbrelations}\,-\ref{eq:syminvrelations})
below.

\section{Inverse braid monoids and inverse mapping class monoids\label{sec:ib_imc}}

\begin{figure}
\epsfbox{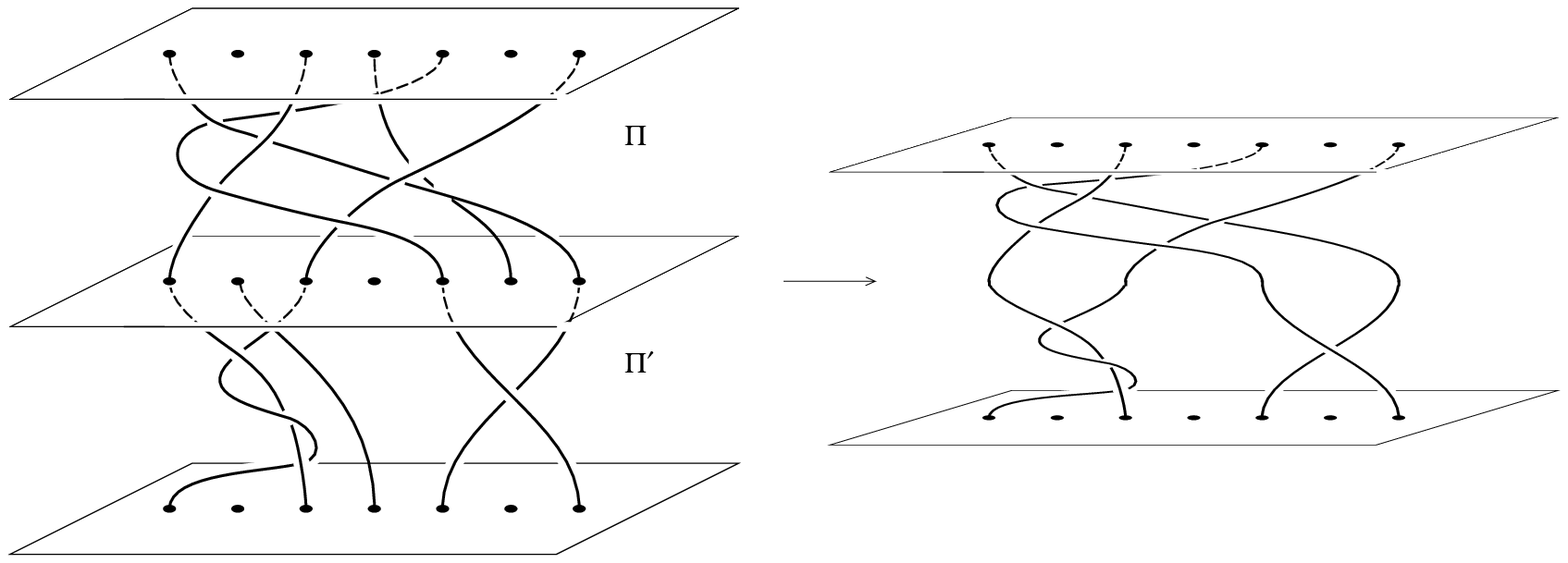}
\caption{} \label{fi:vv2}
\end{figure}

Inverse braid monoid arises from 
a very natural operation on braids: deleting  several strings.
By the application of this procedure to braids in
$Br_n$ we get {\it partial} braids  \cite{EL}.
Applying the standard procedure of multiplication of braids (concatenation) 
to partial braids we need to eliminate all not complete strings, as it
 is shown at Figure~\ref{fi:vv2}. If all strings are eliminated we get
 an {\it empty} braid.  The set of all
partial braids  with this operation forms an inverse braid monoid $IB_n$. 

Usually  the braid group $Br_n$ 
is given by the following Artin presentation \cite{Art1}.
It has the generators $\sigma_i$, 
$i=1, ..., n-1$, and two types of relations: 
\begin{equation}
 \begin{cases} \sigma_i \sigma_j &=\sigma_j \, \sigma_i, \ \
\text{if} \ \ |i-j|
>1,
\\ \sigma_i \sigma_{i+1} \sigma_i &= \sigma_{i+1} \sigma_i \sigma_{i+1} \, .
\end{cases} \label{eq:brelations}
\end{equation}

There exist other presentations of the braid group.
Let 
\begin{equation*}
\sigma = \sigma_1 \sigma_{2} \dots \sigma_{n-1},
\label{eq:sigma}
\end{equation*} 
then the group $Br_n$ is generated by $\sigma_1$ and $\sigma$ because
\begin{equation*}
\sigma_{i+1} =\sigma^i \sigma_1 \sigma^{-i}, \quad i =1, \dots
{n-2}.
\label{eq:sigma_i}
\end{equation*} 
The relations for the generators $\sigma_1$ and $\sigma$ are the 
following
\begin{equation}
 \begin{cases}
\sigma_1 \sigma^i \sigma_1 \sigma^{-i} &= 
\sigma^i \sigma_1 \sigma^{-i} \sigma_1 \ \  \text{for} \ \
2 \leq i\leq {n / 2}, \\
\sigma^n &= (\sigma \sigma_1)^{n-1}.
\end{cases} \label{eq:2relations}
\end{equation}
The presentation (\ref{eq:2relations}) was given by Artin in the initial 
paper \cite{Art1}.
This presentation was also mentioned in the books by F.~Klein \cite{Kl}
and by H.~S.~M.~Coxeter and W.~O.~J.~Moser \cite{CM}.

Classical braid group $Br_n$ can be defined also as 
the mapping class group of a disc $D^2$ with $n$ points deleted (or fixed) and 
with its boundary fixed, or as the subgroup of the
automorphism group of a free group 
 $\operatorname{Aut} F_n, $ 
generated by the following automorphisms:
\begin{equation*} 
\begin{cases} 
x_i &\mapsto x_{i+1},
\\ x_{i+1} &\mapsto x_{i+1}^{-1}x_ix_{i+1}, \\
x_j &\mapsto x_j, j\not=i,i+1. 
\end{cases} \label{eq:autf}
\end{equation*}
 Geometrically
this action is depicted in Figure~\ref{fi:mapcl},
where
$x_i$ correspond to the canonical loops on $D^2$ which form the
generators of the fundamental group of the punctured disc.

This topological interpretation of the braid group was continued in 
\cite{Ve11} for the inverse braid monoid. Let $Q_n$ the set of 
$n$ fixed  points of  a disc $D^2$.
The fundamental group of $D^2$ with these points deleted is isomorphic to $F_n$. 
Consider homeomorphisms
 of $D^2$ onto a copy of the same disc with the condition that
only $k$ points of $Q_n$,  $k \leq n$ (say $i_1, \dots, i_k$) are mapped
bijectively onto the $k$ points (say $j_1, \dots, j_k$) of the second copy 
of $D^2$. 
Consider the isotopy 
classes of such homeomorphisms and denote the set of them by $IM_n(D^2)$. Evidently it is a monoid.

\begin{Theorem} {\rm \cite{Ve11}} The monoids $IB_n$ and  $IM_n(D^2)$
are isomorphic.
\end{Theorem} 

\begin{figure}
\input MAPCL.TEX
\caption{}\label{fi:mapcl}
\end{figure}

The following presentation for the inverse braid monoid was obtained in
\cite{EL}. It has the generators $\sigma_i, \sigma_i^{-1} $, $i=1,\dots,n-1,$
$\epsilon$, and relations
\begin{equation}
 \begin{cases} 
&\sigma_i\sigma_i^{-1}=\sigma_i^{-1}\sigma_i =1, \ \text {for \ all} \ i, \\
&\epsilon \sigma_i  =\, \sigma_i \epsilon \ \ \text {for } i\geq 2,   \\ 
&\epsilon\sigma_1 \epsilon  = \sigma_{1} \epsilon \sigma_1 \epsilon = 
\epsilon\sigma_{1} \epsilon \sigma_1, \\
&\epsilon = \epsilon^2 = \epsilon \sigma_1^2= \sigma_1^2 \epsilon
\end{cases} \label{eq:invbrelations}
\end{equation}
and the braid relations (\ref{eq:brelations}).

Geometrically the generator $\epsilon$ means that the first string in the trivial braid is absent.
 
If we replace the first relation in (\ref{eq:invbrelations})
 by the following set of relations
\begin{equation}
\sigma_i^2 =1, \ \text {for \ all} \ i, \\
 \label{eq:syminvrelations}
\end{equation}
and delete the superfluous relations 
$$\epsilon =  \epsilon \sigma_1^2= \sigma_1^2 \epsilon, $$
we get a presentation of the symmetric inverse monoid $I_n$ \cite{Po} . 
We also can simply add the relations (\ref{eq:syminvrelations})
if we do not worry about redundant relations.
We get a canonical map \cite{EL}
\begin{equation}
\tau_n: IB_n\to I_n,
 \label{eq:tauIBn}
\end{equation}
which is a natural extension of the corresponding map for the braid 
and symmetric groups.

More balanced relations for the inverse braid monoid were 
obtained in \cite{Gil}.
Let $\epsilon_i$ denote the trivial braid with $i$th string deleted,
formally:
\begin{equation*}
\begin{cases} \epsilon_1 &= \epsilon, \\
\epsilon_{i+1} &= \sigma_i^{\pm 1}\epsilon_{i}\sigma_i^{\pm 1}.  
\end{cases} \end{equation*} 
The generators are: $\sigma_i, \sigma_i^{-1} $, $i=1,\dots,n-1,$
$\epsilon_i$, $i=1,\dots,n$, and relations  are the following:

\begin{equation}
 \begin{cases} 
&\sigma_i\sigma_i^{-1}=\sigma_i^{-1}\sigma_i =1, \ \text {for all} \ i, \\
&\epsilon_j \sigma_i =\, \sigma_i \epsilon_j \ \ \text {for } \ j \not= i, i+1, \\ 
&\epsilon_i\sigma_i =  \sigma_{i} \epsilon_{i+1},  \\
&\epsilon_{i+1}\sigma_i =  \sigma_{i} \epsilon_{i},  \\
&\epsilon_i = \epsilon_i^2 , \\
& \epsilon_{i+1} \sigma_i^2= \sigma_i^2 \epsilon_{i+1} = \epsilon_{i+1}, \\
&\epsilon_i \epsilon_{i+1} \sigma_i = \sigma_{i} \epsilon_i \epsilon_{i+1}
=\epsilon_i\epsilon_{i+1},
\end{cases} \label{eq:invbrelations2}
\end{equation}
plus the braid relations (\ref{eq:brelations}).

If we take the 
presentation  (\ref{eq:2relations}) for the braid group we get a presentation
of the inverse braid monoid with generators $\sigma_1, \sigma_1^{-1}, \sigma, \sigma^{-1}$,  $\epsilon$,
and relations \cite{Ve11}:
\begin{equation}
 \begin{cases} 
&\sigma_1\sigma_1^{-1}=\sigma_1^{-1}\sigma_1 =1,  \\
&\sigma\sigma^{-1}=\sigma^{-1}\sigma =1,  \\
&\epsilon \sigma^{i}\sigma_1\sigma^{-i}  
=\, \sigma^{i}\sigma_1\sigma^{-i} \epsilon \ \ \text {for } 1\leq i \leq n-2,   \\ 
&\epsilon\sigma_1 \epsilon  = \sigma_{1} \epsilon \sigma_1 \epsilon = 
\epsilon\sigma_{1} \epsilon \sigma_1, \\
&\epsilon = \epsilon^2 = \epsilon \sigma_1^2= \sigma_1^2 \epsilon,
\end{cases} \label{eq:3invbrelations}
\end{equation}
plus (\ref{eq:2relations}).

Similar to  braids the notion of {\it mapping class monoid}
was introduced in \cite{Ve11}. The definition is as follows.

It is convenient to consider the surface $S_{g,b,n}$, sometimes
with $n$ points fixed, sometimes deleted.
Let $f$ be a homeomorphism of $S_{g,b,n}$ 
which maps  $k$ points,  $k \leq n$, from $Q_n$:  $\{i_1, \dots, i_k\}$ to $k$ points 
$\{j_1, \dots, j_k\}$ also from $Q_n$. The same way let $h$  
be a homeomorphism of $S_{g,b,n}$ 
which maps  $l$ points, $l\leq n$, from $Q_n$, say $\{s_1, \dots, s_l\}$ to $l$ points 
$\{t_1, \dots, t_l\}$ again from $Q_n$. Consider the intersection of the sets
$\{j_1, \dots, j_k\}$ and $\{s_1, \dots, s_l\}$,  let it be the set of cardinality $m$,
it may be empty. Then the composition of $f$ and $h$ maps $m$ points of $Q_n$
to $m$ points (may be different) of $Q_n$. If $m=0$ then the composition have
no relation to the set $Q_n$. Denote the set of isotopy classes of such maps
by $\mathcal I \mathcal {M}_{g,b,n}$. Composition defines a structure of monoid
on $\mathcal I \mathcal {M}_{g,b,n}$.
It is evident  that the monoid $\mathcal I \mathcal {M}_{g,b,n}$ is inverse, so we call it the {\it inverse mapping class 
monoid}. 
If $g=0 $ and $b=1$ we get the inverse braid monoid. In 
the general case of $\mathcal{ IM}_{g,b,n}$ the role of the empty braid 
plays the
mapping class group  $ \mathcal {M}_{g,b}$ (without fixed points). 
Each element $f \in \mathcal{ IM}_{g,b,n}$ corresponds a (partial)
bijection 
\begin{equation*}
\{i_1, \dots, i_k\} \to \{j_1, \dots, j_k\}.
\end{equation*}
This defines a canonical homomorphism to the symmetric inverse monoid:
\begin{equation}
\tau_{g,b,n}: \mathcal{ IM}_{g,b,n}\to I_n.
 \label{eq:tauIMn}
\end{equation}

We remind that a monoid $M$ is {\it factorisable} if $M= EG$ where $E $ 
is a set of idempotents of $M$ and $G$ is a subgroup of $M$. 

\begin{Proposition} \cite{Ve11} The monoid
$\mathcal I \mathcal {M}_{g,b,n}$ 
can be written in the form
$$\mathcal I \mathcal {M}_{g,b,n} = E \mathcal {M}_{g,b,n},$$
where $E $ is a set of 
idempotents of $\mathcal I \mathcal {M}_{g,b,n}$  and $\mathcal {M}_{g,b,n}$ is 
the corresponding mapping class group. So, this monoid is factorisable.
\end{Proposition}

Similar to braids we introduce the idempotent elements 
$\epsilon_i \in \mathcal{IM}_{g,b,n}$ as isotopy classes of identity
map 
\begin{equation*}
Id: S_{g,b} \to S_{g,b},
\end{equation*}
where during isotopy all points of the set $Q_n$ are fixed with the 
exception of the point with number $i$. The element $\epsilon_1$
is denoted by $\epsilon$.

We call an element $m$ of the mapping class monoid 
$ \mathcal{IM}_{g,b,n}$ $i$-{\it Makanin} or $i$-{\it Brunnian} if it
satisfies  the equation:
\begin{equation}
\epsilon_i m = \epsilon_i.
\label{eq:imakm}
\end{equation}
Geometrically this means that if we fill the deleted point $i$ at the surface,
then a homeomorphism $h$ lying in the class $m$
\begin{equation*}
h: S_{g,b} \to S_{g,b},
\end{equation*}
becomes isotopical to the identity map
\begin{equation*}
Id: S_{g,b}/\{1, \dots, \hat{i}, \dots, n\} \to 
S_{g,b}/\{1, \dots, \hat{i}, \dots, n\}.
\end{equation*}
 The condition (\ref{eq:imakm}) is equivalent to the condition
\begin{equation*}
m \epsilon_{\tau (m)(i)}  = \epsilon_{\tau(m)(i)},
\label{eq:imakm2}
\end{equation*}
where $\tau$ is the canonical map to the symmetric inverse monoid
 (\ref{eq:tauIMn}).
With the exception of $\epsilon_i$ itself all such elements belong to
the mapping class group $ \mathcal{M}_{g,b,n}$.
 We  denote the subgroup
of $i$\,-Makanin elements of the mapping class group by $A_i$. The subgroups $A_i$, $i=1, \dots, n$,
 are conjugate. 
The intersection of all subgroups of $i$\,-Makanin elements is the of 
{\it Makanin} or {\it Brunnian} subgroup of the mapping class group
\begin{equation*} 
Mak_{g,b,n} =\cap_ {i=1}^{n} A_i.
\label{eq:mak}
\end{equation*}
That is the same as $m\in Mak_{g,b,n}$ if and only if the equation (\ref{eq:imakm})
holds for all $i$.
Certain properties of {\it Makanin}  subgroups of the mapping class groups
are studied in the work Benson Farb, Christopher J. Leininger and 
Dan Margalit \cite{FLM}

The purpose of this paper is to develop further the theory of
 inverse mapping class monoids, to 
demonstrate that canonical properties of 
 mapping class  groups 
 %and notions based on braids 
 often have there smooth continuation 
for the inverse mapping class 
monoid $\mathcal{IM}_{g,b,n}$.

\section{Properties of inverse mapping class monoids\label{sec:propm}}

Following the ideology of considering partial symmetries instead of global
ones \cite{Law} we can define the monoid of partial automorphisms of a 
group as follows. For a group $G$ consider a set of partial isomorphisms 
\begin{equation*}
f: H\to K,
\end{equation*}
where $H,K$ are subgroups of $G$. The composition of $f$ with the 
isomorphism 
\begin{equation*}
g: L\to M,
\end{equation*}
($L,M$ are subgroups of $G$) 
is a superposition of $f$ and $g$ which is defined on $(K\cap L)f^{-1}$:
 \begin{equation*}
fg: (K\cap L)f^{-1}\to (K\cap L)g,
\end{equation*} 
The set of all such partial isomorphisms with this operation form a
monoid which is evidently inverse and which we call the {\it inverse
partial automorphism monoid of a group} $G$ and denote by 
$IPA(G)$.

In the case when the group $G$ is a finitely generated free group $F_n$
the following submonoid $E F_n$ of $IPA(F_n)$ was defined in 
\cite{Ve11}.
Let $a$ be an element of the symmetric
inverse monoid $I_n$, $a\in I_n$, $J_k =\{j_1, \dots, j_k\}$ 
is the image of $a$, and the elements $i_1, \dots, i_k$ belong to the
domain of the definition of $a$. The monoid $E F_n$
consists of isomorphisms
$$<x_{i_1}, \dots, x_{i_k}> \, \to \ <x_{j_1}, \dots, x_{j_k}>$$
expressed by 
$$f_a :x_i\mapsto w_i^{-1} x_{a(i)}w_i, $$
if $i$ is among $i_1, \dots, i_k$ and not defined otherwise and 
$w_i$ is a word on $x_{j_1}, \dots, x_{j_k}$.
The composition of $f_a$ and $g_b$, $a, b\in I_n$, 
is defined for $x_i$ belonging to the domain of $a\circ b$.
We put $x_{j_m}=1$ in a word $w_i$ if $x_{j_m}$ does not belong
to the domain of definition of $g$.
 We define a map $\phi_n$ from $IB_n$ to $E F_n$
expanding the canonical inclusion 
\begin{equation*}
Br_n \to  \operatorname{Aut} F_n
 \end{equation*}
by the condition that $\phi_n(\epsilon)$ 
as a partial isomorphism of $F_n$ is given by the
formula
\begin{equation*} 
\phi_n(\epsilon)(x_i) = \begin{cases}
x_i {\text{ if} } \ i\geq 2 , \\ 
{\text {not defined,  if }} i=1 . 
\end{cases} \label{eq:endf1}
\end{equation*}

Using the presentation (\ref{eq:invbrelations}) we see that $\phi_n$ is 
correctly defined homomorphism of monoids
\begin{equation*}
\phi_n: IB_n \to  E F_n.
 \end{equation*}
 
The following statement was proved in \cite{Ve11}. 
\begin{Theorem}  The homomorphism $\phi_n$ is a monomorphism.
\label{theo:isoendo} 
\end{Theorem} 

As usual, we define  the $(g, b,n)$-{\it surface group} as a group with 
the presentation
\begin{equation*}
\pi_{g,b,n}= <a_1, c_1, \dots, a_g, c_g, v_1, \dots, v_b, u_1, \dots, u_n \
 | \ \prod_{i=1}^{n}u_i \prod_{l=1}^{b}v_l \prod_{m=1}^{g}[a_m,c_m]>.  
 \end{equation*}
 
In our construction of the inverse mapping class monoid 
$\mathcal{ IM}_{g,b,n}$ we are passing from maps which are fixing 
a set of $n$
points to maps which are fixing the set with a smaller number of points. 
If we consider
these points deleted it means we are filling the holes. On the level of
 fundamental groups this means passing to factor group.
 
Let $H$ be a factor group of $\pi_{g,b,n}$, defined  by the conditions
\begin{equation*}
u_i=1 \text{ for all } \ i\not\in \{i_1, \dots, i_k\},
\end{equation*}
and let $K$ be a factor group of $\pi_{g,b,n}$, defined  by the conditions
\begin{equation*}
u_j=1 \text{  for all } \ j\not\in \{j_1, \dots, j_k\}.
\end{equation*}
Let $t$ be an element of the symmetric
inverse monoid $I_n$, $t\in I_n$, $J_k =\{j_1, \dots, j_k\}$ 
is the image of $t$, and elements $i_1, \dots, i_k$ belong to the
domain of the definition of $t$. Let $w_i$ be a word on letters $a_1, c_1, \dots, a_g, c_g, v_1, \dots, v_b, u_{j_1}, \dots, u_{j_k}$.
Define $I Aut \, \pi_{g,b,n}$ consisting of isomorphisms
\begin{equation*}
f_t: H\to K,
\end{equation*}
such that 
\begin{equation*}
\begin{cases} 
f_t(v_m)= v_m \ \text{for} \ m=1, \dots, b,\\
f_t(u_i) = w_i^{-1} u_{(i)t}w_i, \ \text{if} \ i \ \text{is among} \ i_1, \dots, i_k, \\
%  f_t(u_i)  \ \text{not defined otherwise}.
\end{cases} 
 \end{equation*}
 
 \smallskip\noindent 
for all subgroups $H$ and $K$ of the type defined above, for all 
$k= 0, 1, \dots, n$. 
%Let $In\pi_{g,b,n}$ be the group of inner automorphisms

The composition of $f_t$ and $g_s$, $t, s\in I_n$, 
is defined for $u_i$ belonging to the domain of $t\circ s$.
We put $u_{j_m}=1$ in a word $w_i$ if $u_{j_m}$ does not belong
to the domain of definition of $g_s$.

By our construction to each element of $f_t \in I Aut \, \pi_{g,b,n}$
the element $t\in I_n$ is associated. This gives a canonical 
homomorphism 
\begin{equation}
\tau_{g,b,n}: I Aut \, \pi_{g,b,n}\to I_n.
 \label{eq:tauIAn}
\end{equation}

For the simplicity of exposition we restrict ourselves to the case
of empty boundary which was described in the book of W.~Magnus,
 A.~Karrass and D.~Solitar \cite{MKS}, Sec.~3.7, and we denote
$\pi_{g,0,n}$ by $\pi_{g,n}$.
Let us define an equivalence relation in $I Aut \, \pi_{g,n}$: for 
$f_1$ and $f_2$ 
\begin{equation*}
f_1, f_2: H\to K.
\end{equation*}
We are considering the homeomorphism of a surface with $k$ points
deleted (among our $n$ fixed points) onto another copy of the surface
with (may be) different $k$ points deleted 
\begin{equation*}
h: S_{g,n}/ \{i_1, \dots, i_k\} \to S_{g,n}/ \{j_1, \dots, j_k\}.
\end{equation*}
In such a case
it is reasonable to consider the bijection between the conjugacy classes
of $\pi_1(S_{g,n}/ \{i_1, \dots, i_k\})$ and $\pi_1(S_{g,n}/ \{j_1, \dots, j_k\})$ as an algebraic image of $h$.
This equivalence is in fact a congruence on the monoid
$I Aut \, \pi_{g,n}$ and we denote the corresponding factor monoid by
$IOut \, \pi_{g,n}$. This is similar to the classical case of groups when
it is necessary to factorize by  inner automorphisms. 
The homomorphism $\tau_{g,b,n}: I Aut \, \pi_{g,b,n}\to I_n $ of
(\ref{eq:tauIAn}) factors through $IOut \, \pi_{g,n}$ and we have a 
homomorphism
\begin{equation}
\tau_{g,b,n}: I Out \, \pi_{g,b,n}\to I_n.
 \label{eq:tauIOn}
\end{equation}

These considerations also define a homomorphism of monoids
\begin{equation*}
\psi_{g,n}: \mathcal I \mathcal {M}_{g,n} \to  IOut \, \pi_{g,n}.
 \end{equation*}
which is compatible with homomorphisms $\tau_{g,b,n}$ of
(\ref{eq:tauIMn}) and (\ref{eq:tauIOn}).
\begin{Theorem}  The homomorphism $\psi_n$ is an isomorphism of monoids.
\label{theo:iso_mcm} 
\end{Theorem} 
\begin{proof}
Monoid $\mathcal{ IM}_{g,n}$ as a set is a disjoint union of copies of 
mapping class groups
$\mathcal{M}_{g,k}$, for $k= 0, \dots n$, namely, $C^k_n$ copies of $\mathcal{M}_{g,k}$ for each $k$. 
Our constructions are done so that for each copy of $\mathcal{M}_{g,k}$ 
the map $\psi_{g,n}$ is a bijection because of the Dehn-Nilsen-Baer theorem. 
\end{proof}

\section{Inverse mapping class monoids for punctured sphere\label{sec:sph}}

The following presentation for the 
mapping class group of a punctured sphere $\mathcal {M}_{0,n}$
was obtained by W.~Magnus \cite{Mag1}, see also\cite{MKS}. 
Let $\sigma_1, \dots, \sigma_{n-1}$, denote the classes of 
homeomorphisms such that $\sigma_i$ locally interchanges the points 
with numbers $i$ and $i+1$.  Then the presentation has 
generators $\sigma_1, \dots, \sigma_{n-1}$, which satisfy the braid 
relations (\ref{eq:brelations}), the sphere relation

\begin{equation}
\sigma_1 \sigma_2 \dots \sigma_{n-2}\sigma_{n-1}^2\sigma_{n-2} \dots
\sigma_2\sigma_1 =1
\label{eq:spherelation}
\end{equation}
and the following relation

\begin{equation}
(\sigma_1 \sigma_2 \dots \sigma_{n-2}\sigma_{n-1})^n =1.
\label{eq:sphe_mc}
\end{equation}
Let $\Delta$ be the Garside's   fundamental word   in the braid 
group $Br_{n}$ \cite{Gar}. It can be in particular defined by the formula:
$$\Delta = \sigma_1 \dots \sigma_{n-1} \sigma_1 \dots \sigma_{n-2} \dots  
\sigma_1 \sigma_2 \sigma_1.$$
If we use Garside's notation $\Pi_t\equiv \sigma_1\dots \sigma_t$, then
$\Delta \equiv \Pi_{n-1} \dots \Pi_1$.
If the generators $\sigma_1,$ $\sigma_2$, $\dots$, $\sigma_{n-2},$
$\sigma_{n-1}$,
are subject to the braid relations (\ref{eq:brelations}), then the 
condition (\ref{eq:sphe_mc}) is equivalent to the following
\begin{equation*}
\Delta^2 =1.
\label{eq:sphe_mcD}
\end{equation*}

If we consider the presentation of the braid group with two generators the
sphere relation has the form 

\begin{equation*}
\sigma^n (\sigma_1^{-1}\sigma)^{1-n} =1,
\label{eq:spherelation2g}
\end{equation*}
and the sphere mapping class relation is
\begin{equation*}
\sigma^n  =1.
\label{eq:sphemc2g}
\end{equation*}
So, the mapping class group $\mathcal {M}_{0,n}$ in two generators is described 
by the braid relations   and the following two relations
\begin{equation}
\begin{cases}
\sigma^n  =1, \\
(\sigma_1^{-1}\sigma)^{n-1} =1.
\end{cases}
\label{eq:sphe_mc2r}
\end{equation}

\begin{Theorem} We get a presentation of the inverse mapping class 
monoid for punctured sphere $\mathcal{IM}_{0,n}$ if we take the generators
$\sigma_1,$ $\sigma_2$, $\dots$, $\sigma_{n-2},$
$\sigma_{n-1}$, $\sigma_1^{-1},$ $\sigma_2^{-1}$, $\dots$, $\sigma_{n-2}^{-1},$
$\sigma_{n-1}^{-1}$, $\epsilon$ (or $\epsilon_1,$ $\epsilon_2$, $\dots$, $\epsilon_{n-1},$ $\epsilon_{n}$ instead of one $\epsilon$) subject to 
the sphere braid relations (\ref{eq:brelations}, (\ref{eq:spherelation}),
the sphere mapping class relation (\ref{eq:sphe_mc}) 
and the inverse braid relations (\ref{eq:invbrelations}) (or 
 (\ref{eq:invbrelations2}).
 \label{theo:presimcm}
\end{Theorem} 
\begin{proof} 
We use the fact that the mapping class group $\mathcal{M}_{0,0,n}$ is a 
factor  group of the $n$-string braid group for sphere and the ideas of 
the proof of the presentation (\ref{eq:invbrelations}) for the inverse 
braid monoid $IB_n$ in the work  \cite{EL}. 

Denote temporarily by $P_n$ the monoid defined by the presentation
by the generators  
$\sigma_1,$ $\sigma_2$, $\dots$, $\sigma_{n-2},$
$\sigma_{n-1}$, $\sigma_1^{-1},$ $\sigma_2^{-1}$, $\dots$, $\sigma_{n-2}^{-1},$
$\sigma_{n-1}^{-1}$, $\epsilon$ 
the sphere braid relations (\ref{eq:brelations}, (\ref{eq:spherelation}),
the sphere mapping class relation (\ref{eq:sphe_mc}) 
and the inverse braid relations (\ref{eq:invbrelations}). To define a map 
$$\Theta : P_n\to \mathcal{M}_{0,n}$$
we associate to each word in the alphabet on letters
$\sigma_1,$ $\sigma_2$, $\dots$, $\sigma_{n-2},$
$\sigma_{n-1}$, $\sigma_1^{-1},$ $\sigma_2^{-1}$, $\dots$, $\sigma_{n-2}^{-1},$
$\sigma_{n-1}^{-1}$, $\epsilon$ the corresponding composition of 
(classes of) homeomorphisms. To juxtaposition of words there evidently
corresponds the composition of (classes of) homeomorphisms. 
The fact that $\Theta $ is well defined 
(respects the relations) for some relations are classical facts, for others
it is evident, probably, with the exception of the third relation in  
(\ref{eq:invbrelations}). The equalities of this relation follow from
the fact that all four given classes (containing identity) ignore the first 
two points. 

Let $\epsilon_{k+1, n}$ denote the the isotopy class of of a homeomorphism
that fixes the first $k$ point and does not care about the rest $n-k$
points. On the level of braids (on a sphere) it corresponds to the
partial braid with the trivial first $k$ strings 
and the absent rest $n-k$ strings. The element $\epsilon_{k+1, n}$ 
 can be expressed using the generator $\epsilon$ 
or the generators $\epsilon_i$ as follows
\begin{equation} \epsilon_{k+1, n}=
\epsilon\sigma_{n-1}\dots\sigma_{k+1}\epsilon \sigma_{n-1}\dots\sigma_{k+2}
\epsilon\dots \epsilon\sigma_{n-1}\sigma_{n-2}\epsilon \sigma_{n-1}\epsilon,
 %\label{eq:form_inv}
\end{equation}
\begin{equation*} \epsilon_{k+1, n}=
\epsilon_{k+1}\epsilon_{k+2} \dots \epsilon_{n},
 \label{eq:espki}
\end{equation*}
From our construction of the inverse mapping class monoid it 
follows that every element of $\mathcal{IM}_{0,n}$ 
represented 
by a homeomorphism $h$ of $S_{g,b,n}$ which maps  $k$ points,  $k \leq n$, from $Q_n$:  $\{i_1, \dots, i_k\}$ to $k$ points 
$\{j_1, \dots, j_k\}$ from $Q_n$
can be expressed in
the form
\begin{equation} \sigma_{i_1}\dots\sigma_{1}\dots \sigma_{i_k}\dots
\sigma_{k} \,  \epsilon_{k+1, n} \, x \,  \epsilon_{k+1,n} \, \sigma_{k}\dots\sigma_{j_k}\dots\sigma_{1}
\dots\sigma_{j_1}, \\ 
\label{eq:form_inv}\end{equation}
\begin{equation*} k\in \{0,\dots, n\}, \ x\in \mathcal{M}_{0,k}, \ 
0\leq i_1<\dots<i_k\leq n-1  \ 
\text{and} \  0\leq j_1<\dots<j_k\leq n-1.
\end{equation*}
Geometrically this means that we first send the points $\{i_1$, $\dots$,
$i_k\} $ to the points $ \{1, \dots, k\}$, then apply a homeomorphism from
the mapping class {\it group}  $\mathcal{M}_{0,k}$, and then 
send the points $\{1$, $\dots$, $k\} $ to the points 
$ \{j_1, \dots, j_k\}$.

Note that in the formula (\ref{eq:form_inv}) we can  delete one of the 
$\epsilon_{k+1,n}$, but we shall use the form (\ref{eq:form_inv}) because of 
convenience: two symbols $\epsilon_{k+1,n}$ serve as markers to distinguish
the elements of $\mathcal{M}_{0,k}$.

The element $x$ belongs to the mapping class group $\mathcal{M}_{0,k}$, so 
it can be expressed as a word on the letters
$\sigma_1,$ $\sigma_2$, $\dots$, 
$\sigma_{k-1}$, $\sigma_1^{-1},$ $\sigma_2^{-1}$, $\dots$, 
$\sigma_{k-1}^{-1}$. This proves that homomorphism $\Theta$ is onto.

Let us prove that $\Theta$ is a monomorphism. Suppose that for two words
$W_1, W_2 \in P_n$ we have
$$\Theta(W_1) = \Theta (W_2).$$
It means that the corresponding homeomorphisms both map the  set 
of points $\{i_1, \dots, i_k\}$ onto the set of points $\{j_1, \dots, j_k\}$ 
and they are isotopic in the class of homeomorphisms mapping
$\{i_1, \dots, i_k\}$ onto  $\{j_1, \dots, j_k\}$.
Using relations (\ref{eq:brelations}) and
(\ref{eq:invbrelations}) the same way as in \cite{EL} transform the words
$W_1, W_2$ into the form (\ref{eq:form_inv}) 
$$\sigma(i_1,\dots i_k; k) \epsilon_{k+1,n}x \epsilon_{k+1,n}
\sigma(k; j_1,\dots j_k).$$ 
Then  the corresponding fragments $\sigma(i_1,\dots i_k; k)$ and
$\sigma(k, j_1,\dots j_k; k)$ for $W_1$ and $W_2$ coincide.
The  elements $x_1$ of $W_1$ and $x_2 $ of $W_2$, which are the words 
on $\sigma_1, \dots, \sigma_k$ and $\sigma_1^{-1}, \dots, \sigma_k^{-1}$,
 correspond after $\Theta$ to isotopic homeomorphisms 
\begin{equation*}
h_1, h_2: S_{0,k} \to S_{0,k},
\end{equation*} 
Hence   $x_1$ can be transformed into $x_2$ using  relations for the 
  mapping class group $\mathcal{M}_{0,k}$. So, the words $W_1$ and $W_2$
represent the same element in $P_n$. 
\end{proof}

\begin{Proposition} We get a presentation of the inverse mapping class 
monoid for punctured sphere $\mathcal{IM}_{0,n}$ if we take the generators
$\sigma_1,$ $\sigma$, 
 $\sigma_1^{-1},$ $\sigma^{-1}$,  $\epsilon$  subject to 
the  relations (\ref{eq:2relations}),
 (\ref{eq:3invbrelations})
and   (\ref{eq:sphe_mc2r}) .
\end{Proposition} 
\hfill $\square$

%\begin{Proposition} 
The generators $\epsilon_i$ commute with $\Delta$ in
the following way \cite{Ve11}:
\begin{equation*}
\epsilon_i\Delta = \Delta \epsilon_{n+1-i}.
\end{equation*}
%\end{Proposition}
%\hfill $\square$

Let $\mathcal E$ be the monoid generated by one idempotent generator $\epsilon$ . 
\begin{Proposition} The abelianization of  $\mathcal{M}_{0,n}$, $n\geq 2$, is isomorphic to
$\mathcal E \oplus \Z/m\Z$, where $m= 2(n-1)$ i $n$ is even and $m= (n-1)$
if $n$ is odd. The canonical map of abelianization
\begin{equation*}
a: \mathcal{M}_{0,n} \to \mathcal E \oplus \Z/m\Z
\end{equation*}
is given by the formula:
\begin{equation*}
\begin{cases}
a(\epsilon_i) = \epsilon ,\\
a(\sigma_i) = \bar{ 1} .
\end{cases}
\end{equation*}
\end{Proposition}
\hfill $\square$
 
The mapping class groups $\mathcal{M}_{0,n}$ are centerless for $n\geq 3$
\cite{GVB}. So this gives the possibility to understand the center of 
$\mathcal{IM}_{0,n}$. Denote by $\mathcal{M}_{0, \, n ; \, i; \, j}$  
(one-element) subsets of $\mathcal{IM}_{0,n}$ consisting of the class
of homeomorphisms with the only condition that they send the point $i$
to the point $j$ and denote by 
$\mathcal{M}_{0, \,n ; \, i_1, i_2; \, j_1, j_2}$  
(two-element) subsets of $\mathcal{IM}_{0,n}$ consisting of the class
of homeomorphisms with the only condition that they send the two-point 
set $\{i_1, i_2\}$ to the two-point set $\{j_1, j_2\}$ (cf. Examples 2 ans 3 
below).

\begin{Proposition} The center of $\mathcal{IM}_{0,n}$ consists of the 
disjoin union 
%\begin{equation}
$$
C(\mathcal{IM}_{0,n}) = \mathcal{M}_{0,0} \amalg (\amalg_{i,j} 
\mathcal{M}_{0, \, n ; \, i; \, j}) \amalg 
(\amalg_{\{i_1, i_2\}, \{j_1, j_2\}}  
\mathcal{M}_{0, \, n ; \, i_1, i_2; \, j_1, j_2}
$$
%\end{equation}
\end{Proposition}
\hfill $\square$

\section{The word problem\label{sec:wordprob}}

R.~Gillette and J.~Van Buskirk in \cite{GVB} obtained an analogue of the Markov 
normal form for the sphere braid groups $Br_n(S^2)$ and for the punctured 
sphere mapping class group $\mathcal{M}_{0, \, n}$. This form can be obtained
algorithmically, so, it gives a solution of the word problem. These ideas
can be applied to the cases of inverse braid and mapping class monoids.

Let us remind the main points of the Markov's construction. 
Define the elements $s_{i,j}$, $1\leq i<j\leq m$, of the classical
braid group $Br_m$ by the formula:
$$s_{i,j}=\sigma_{j-1}...\sigma_{i+1}\sigma_{i}^2\sigma_{i+1}^{-1}...
\sigma_{j-1}^{-1}.$$ 
These elements satisfy the following Burau
relations: 
\begin{multline}
\begin{cases}s_{i,j}s_{k,l}=s_{k,l}s_{i,j}
\ \text {for} \ i<j<k<l \ \text {and} \ i<k<l<j, \\
s_{i,j}s_{i,k}s_{j,k}=s_{i,k}s_{j,k}s_{i,j} \ \text {for} \
i<j<k, \\
s_{i,k}s_{j,k}s_{i,j}=s_{j,k}s_{i,j}s_{i,k} \ \text
{for} \ i<j<k, \\
s_{i,k}s_{j,k}s_{j,l}s_{j,k}^{-1}=s_{j,k}s_{j,l}s_{j,k}^{-1}s_{i,k}
\ \text {for} \ i<j<k<l.\\
\end{cases}
\label{eq:burau}
\end{multline}
 The elements
$s_{i,j}$ with the relations (\ref{eq:burau}) give a~presentation of the pure
braid group $P_m$ \cite{Mar2}.

 Let us define the elements
$\sigma_{k,l}$, $1\leq k\leq l\leq m$ by the formulas
$$\sigma_{k,k}=e,$$ 
$$\sigma_{k,l}=\sigma_{k}^{-1}...\sigma_{l-1}^{-1}.$$
We denote by the same symbols the images of the defined elements in the
braid group of a sphere. 
Let $w(x_1,..., x_m$) be a word with (possibly empty) entries of
$x_i^\delta$, where $x_i$ are some letters and $\delta$ may be
$\pm 1$. 
\begin{Theorem} \cite{GVB}
 Every element of the group $Br_{n}(S^2)$ can be uniquely written in
the form
\begin{equation*}
\sigma_{i_{n},n}...\sigma_{i_{j},j}...\sigma_{i_{2},2} w_{1}(s_{1,2},...,s_{1,n-1})...w_j(s_{j,j+1},...,s_{j,n-1})...
w_{n-3}(s_{n-3,n-2}, s_{n-2,n-1})\Delta^{2\delta}, \label{eq:markforb}
\end{equation*}
where elements $s_{j,j+1},...,s_{j,n-1}$ generate a free group
and $\delta$ may be $\pm 1$. There exists an algorithm of obtaining 
the normal form for any word on letters $\sigma_1, \dots, \sigma_{n-1}$,
which gives a solution of the word problem for $Br_{n}(S^2)$. 
\label{Theorem:mnfb}
\end{Theorem}
\begin{Theorem} \cite{GVB}
 Every element of the group $\mathcal{M}_{0,n}$  can be uniquely written in
the form
\begin{equation*}
\sigma_{i_{n},n}...\sigma_{i_{j},j}...\sigma_{i_{2},2} w_{1}(s_{1,2},...,s_{1,n-1})...w_j(s_{j,j+1},...,s_{j,n-1})...
w_{n-3}(s_{n-3,n-2}, s_{n-2,n-1}), \label{eq:markformc}
\end{equation*}
where elements $s_{j,j+1},...,s_{j,n-1}$ generate a free group
and $\delta$ may be $\pm 1$. There exists an algorithm of obtaining 
the normal form for any word on letters $\sigma_1, \dots, \sigma_{n-1}$,
which gives a solution of the word problem for $\mathcal{M}_{0,n}$. 
\label{Theorem:mnfmc}
\end{Theorem}

\begin{Theorem}
 Every word $W$ in $IBr_{n}(S^2)$ can be uniquely written in
the form 
\begin{equation*} \sigma_{i_1}\dots\sigma_{1}\dots \sigma_{i_k}\dots\sigma_{k}
\epsilon_{k+1, n}x \epsilon_{k+1, n}\sigma_{k}\dots\sigma_{j_k}\dots\sigma_{1}
\dots\sigma_{j_1},
 \\ 
%\label{eq:form_inv}
\end{equation*}
\begin{equation*} k\in \{0,\dots, n\}, x\in Br_k(S^2),
0\leq i_1<\dots<i_k\leq n-1  \ 
\text{and} \  0\leq j_1<\dots<j_k\leq n-1.
\end{equation*} 
where $x$ is written in the Markov normal form for $Br_k(S^2)$
\label{Theorem:mnfib}
\end{Theorem}
\begin{proof}
It follows from the proof of Theorem~3.1 in \cite{Ve11}.
\end{proof}

\begin{Theorem}
 Every word $W$ in $\mathcal{IM}_{0, \, n}$ can be uniquely written in
the form 
\begin{equation*} \sigma_{i_1}\dots\sigma_{1}\dots \sigma_{i_k}\dots\sigma_{k}
\epsilon_{k+1, n}x \epsilon_{k+1, n}\sigma_{k}\dots\sigma_{j_k}\dots\sigma_{1}
\dots\sigma_{j_1},
 \\ 
%\label{eq:form_inv}
\end{equation*}
\begin{equation*} k\in \{0,\dots, n\}, x\in \mathcal{M}_{0, \, n},
0\leq i_1<\dots<i_k\leq n-1  \ 
\text{and} \  0\leq j_1<\dots<j_k\leq n-1.
\end{equation*} 
where $x$ is written in the Markov normal form for $\mathcal{IM}_{0, \, n}$
\label{Theorem:mnfim}
\end{Theorem}
\begin{proof}
It follows from the proof of Theorem \ref{theo:presimcm}.
\end{proof}

\section{Examples\label{sec:ex}}

1. The monoid $\mathcal{IM}_{0,0}$ consists of one element of the identical 
map of a sphere
\begin{equation*}
Id: S_{0,0} \to S_{0,0},
\end{equation*} 

2. The monoid $\mathcal{IM}_{0,1}$ consists of two elements, say, unit $1$ 
and the idempotent $\epsilon$.

3. The mapping class group $\mathcal{M}_{0,2}$ is isomorphic to the 
symmetric group
$\Sigma_2$, so the inverse mapping class  monoid $\mathcal{IM}_{0,2}$
is isomorphic to the symmetric inverse monoid $I_2$.

4. The mapping class group $\mathcal{M}_{0,3}$ is isomorphic to the 
symmetric group
$\Sigma_3$, so the inverse mapping class  monoid $\mathcal{IM}_{0,3}$ is isomorphic to the symmetric inverse monoid $I_3$.

5. The mapping class group of a once-punctured torus $\mathcal{M}_{1,1}$ 
is isomorphic to the three-strand braid group $B_3$ and mapping class group 
of a torus (without punctures) is isomorphic to the unimodular group
$SL_2\Z$ which is considered as a factor group of $B_3$
\begin{equation*}
\rho: B_3 \to SL_2\Z.
\end{equation*}
by the relation $\Delta^4=1$.
 So  the inverse 
mapping class  monoid $\mathcal{IM}_{1,1}$
as a set is a disjoint union 
\begin{equation*}
\mathcal{IM}_{1,1} = B_3 \amalg SL_2\Z.
\end{equation*}
Multiplication of elements $b\in B_3$ and $a\in  SL_2\Z$ is defined as 
follows: 
\begin{equation*}
a\cdot b = a\rho(b), \ b\cdot a= \rho(b) a.
\end{equation*}

\end{document}